\documentclass[runningheads]{llncs}
\usepackage[T1]{fontenc}

\usepackage{color}
\definecolor{keywordcolor}{rgb}{0.7, 0.1, 0.1}
\definecolor{tacticcolor}{rgb}{0.0, 0.1, 0.6}
\definecolor{commentcolor}{rgb}{0.4, 0.4, 0.4}
\definecolor{symbolcolor}{rgb}{0.0, 0.1, 0.6}
\definecolor{sortcolor}{rgb}{0.1, 0.5, 0.1}
\definecolor{attributecolor}{rgb}{0.7, 0.1, 0.1}
\definecolor{backcolour}{rgb}{0.95, 0.95, 0.92}

\usepackage{listings}

\usepackage{amssymb}
\usepackage{amsmath}

\usepackage{graphicx}
\usepackage{wrapfig}
\usepackage{caption}
\usepackage{tikz}
\usetikzlibrary{arrows}

\usepackage{enumitem}

\usepackage{hyperref}
\hypersetup{
    citebordercolor = 1 1 1,
    linkbordercolor = 1 1 1,
    filebordercolor = 1 1 1,
    menubordercolor = 1 1 1,
    urlbordercolor  = 1 1 1,
    colorlinks = true,
    linkcolor  = blue,
    citecolor  = magenta,
    urlcolor   = blue,
}
\usepackage{cleveref}

\newcommand{\R}{\mathbb{R}}
\newcommand{\Z}{\mathbb{Z}}
\newcommand{\C}{\mathbb{C}}

\newcommand{\parenth}[1]{\left( #1 \right)}
\newcommand{\norm}[1]{\left\lVert #1 \right\rVert}
\newcommand{\setst}[2]{\left\{ #1 \; \middle\vert \; #2 \right\}}
\newcommand{\of}[1]{\!\parenth{#1}}
\newcommand{\ssq}{\subseteq}

\newcommand{\Pa}{\mathcal{P}}
\newcommand{\Volof}[1]{\operatorname{Vol}\!\parenth{#1}}
\newcommand{\finite}{\operatorname{finite}}

\newcommand{\sinsq}[1]{\sin^2\!\parenth{#1}}
\newcommand{\diff}[1]{\operatorname{d}\!{#1}}

\renewcommand{\Re}{\operatorname{Re}}
\renewcommand{\Im}{\operatorname{Im}}
\renewcommand{\hat}{\widehat}

\newcommand{\mathlib}{\texttt{mathlib}}
\newcommand{\projectfile}[1]{\href{https://github.com/thefundamentaltheor3m/Sphere-Packing-Lean/blob/main/#1}{\texttt{#1}}}
\newcommand{\PR}[1]{PR \href{https://github.com/thefundamentaltheor3m/Sphere-Packing-Lean/pull/#1}{\#{#1}}}
\newcommand{\sorry}{\lstinline`sorry`}

\newcommand{\labelledpoint}[5]{%
    \node [label={[shift={(#3,#4)}]#5}] at (#1, #2) {$\bullet$};
}

\begin{document}
\title{Progress in Formalizing Sphere Packing in Dimension 8}
\author{Sidharth Hariharan\inst{1}\orcidID{0009-0000-8176-902X} \and
Christopher Birkbeck\inst{2}\orcidID{0000-0002-7546-9028} \and
Seewoo Lee\inst{3}\orcidID{0000-0002-5710-2257} \and
Ho Kiu Gareth Ma\inst{4}\orcidID{0009-0008-2940-8233} \and
Bhavik Mehta\inst{5}\orcidID{0000-0001-7892-7891} \and
Auguste Poiroux\inst{6,7}\orcidID{0009-0005-7868-7185} \and
Maryna Viazovska\inst{7}\orcidID{0000-0001-7567-5254}
}
\authorrunning{S. Hariharan et al.}
\institute{Carnegie Mellon University, Pittsburgh PA 15213, United States of America \and
University of East Anglia, Norwich NR4 7TJ, United Kingdom \and
University of California, Berkeley, Berkeley CA 94720, United States of America \and
University of Warwick, Coventry CV4 7AL, United Kingdom \and
Imperial College London, London SW7 2AZ, United Kingdom \and
Math, Inc, San Francisco CA 94107, United States of America \and
École Polytechnique Fédérale de Lausanne, 1015 Lausanne VD, Switzerland
}
\maketitle
\begin{abstract}
In 2016, Viazovska famously solved the sphere packing problem in dimension $8$, using modular forms to construct a `magic' function satisfying optimality conditions determined by Cohn and Elkies in 2003. In March 2024, Hariharan and Viazovska launched a project to formalize this solution and related mathematical facts in the Lean Theorem Prover. A significant milestone was achieved in February 2026: the result was formally verified, with the final stages of the verification done by Math, Inc.'s autoformalization model `Gauss'. We discuss the techniques used to achieve this milestone, reflect on the unique collaboration between humans and Gauss, and discuss project objectives that remain.

\end{abstract}
\section{Introduction}

The sphere packing problem is a notoriously difficult problem in discrete geometry, asking what the densest arrangement of non-overlapping $n$-spheres is in $\R^n$. When $n=1$, the solution is trivial. When $n=2$, the optimal packing is the $A_2$ lattice packing. The result is often attributed to Thue \cite{Thue}, though many credit Toth \cite{Toth} for providing the first rigorous proof. The story is more interesting when $n=3$, and while the solution was originally conjectured by Kepler \cite{KeplerSnowflake} as far back as 1611, it was not proved till the turn of the century by Tom Hales \cite{HalesKeplerInformal}, who went on to formally verify \cite{HalesKeplerFormal} his heavily computer-assisted, computationally involved proof.

In March 2024, Hariharan and Viazovska launched a project to investigate Viazovska's solution in dimension $8$ \cite{Viazovska8} and follow-up calculations by Lee \cite{Seewoo_Ineq} by formalizing them in the Lean Theorem Prover \cite{Lean_System_Description}, building infrastructure along the way to support ongoing and future work in related areas to the end of creating a readable, maintainable, and reusable proof artifact for this remarkable result.
The work done so far by both human formalizers and the autoformalization model `Gauss' has together led to a complete, correct\footnote{in the sense of being formally verified by the Lean kernel} proof that the optimal sphere packing in $\R^8$ is the $E_8$ lattice packing. This paper discusses the project objectives achieved with this milestone and outlines those that remain.

\section{The Mathematics behind the $8$-Dimensional Sphere Packing Problem}

In this section, we give an overview of the informal argument, listing key definitions and theorem statements, both informally and in Lean. All code in this section was written by human formalizers prior to the Gauss autoformalization.

\subsection{The Road to Stating the Main Theorem}

The primary data making up a sphere packing is a set $X \subseteq \R^d$ of centres and some separation $r > 0$ such that for all $x, y \in X$, $\norm{x - y} \geq r$. We thus \href{https://github.com/thefundamentaltheor3m/Sphere-Packing-Lean/blob/6bd7a08c1fecdfe6d958574082e8474e602ecf76/SpherePacking/Basic/SpherePacking.lean#L37-L41}{defined \lstinline|SpherePacking|} to be the \lstinline|structure| bundling together the following data: the set of centres; the separation; the requirement that the separation is positive (with infrastructure to automatically deduce this where possible); and the fact that the separation does, indeed, separate the centres.
\begin{lstlisting}[language=Lean]
structure SpherePacking (d : ℕ) where
  centers : Set (EuclideanSpace ℝ (Fin d))
  separation : ℝ
  separation_pos : 0 < separation := by positivity
  centers_dist : Pairwise (separation ≤ dist · · : centers → centers → Prop)
\end{lstlisting}
This definition does not actually involve the balls that make up the packing. We thus \href{https://github.com/thefundamentaltheor3m/Sphere-Packing-Lean/blob/6bd7a08c1fecdfe6d958574082e8474e602ecf76/SpherePacking/Basic/SpherePacking.lean#L100-L101}{defined \lstinline|SpherePacking.balls|} to be the quantity $\Pa(X) = \bigcup_{x \in X} B_{d}\of{x, r}$:
\begin{lstlisting}[language=Lean]
abbrev SpherePacking.balls (S : SpherePacking d) : Set (EuclideanSpace ℝ (Fin d)) := ⋃ x : S.centers, Metric.ball x (S.separation / 2)
\end{lstlisting}
The prefix \lstinline`SpherePacking` in the above definition not only places \lstinline`balls` in the \lstinline`SpherePacking` namespace (which is useful for disambiguation) but also allows us to use dot notation, whereby if we wish to refer to the balls of a particular sphere packing \lstinline`S`, we can simply write \lstinline`S.balls`.

Next, we define the \textbf{finite density}\footnote{A more accurate term might be `bounded' density, in the sense of being computed over bounded regions.} of a sphere packing as the proportion of a bounded ball it occupies. We further define the \textbf{density} of a sphere packing to be the $\limsup$ of all finite densities. We denote them $\Delta_{\Pa(X)}^{\finite}$ and $\Delta_{\Pa(X)}$ respectively:
\begin{align*}
    \Delta_{\Pa(X)}^{\finite}(R) := \frac{\Volof{\Pa(X) \cap B_d(0, R)}}{\Volof{B_d(0, R)}}
    \qquad \qquad
    \Delta_{\Pa(X)} := \limsup_{R \to \infty} \Delta_{\Pa(X)}^{\finite}(R)
\end{align*}
We \href{https://github.com/thefundamentaltheor3m/Sphere-Packing-Lean/blob/6bd7a08c1fecdfe6d958574082e8474e602ecf76/SpherePacking/Basic/SpherePacking.lean#L103-L107}{define these in Lean} as follows:
\begin{lstlisting}[language=Lean]
def SpherePacking.finiteDensity (S : SpherePacking d) (R : ℝ) : ℝ≥0∞ :=
  volume (S.balls ∩ ball 0 R) / (volume (ball 0 R))
def SpherePacking.density (S : SpherePacking d) : ℝ≥0∞ := limsup S.finiteDensity atTop
\end{lstlisting}
A \textbf{periodic sphere packing} is a sphere packing that is periodic with respect to the additive action of a lattice on the ambient space (meaning, in particular, that it acts on the set of centres). Our \href{https://github.com/thefundamentaltheor3m/Sphere-Packing-Lean/blob/6bd7a08c1fecdfe6d958574082e8474e602ecf76/SpherePacking/Basic/SpherePacking.lean#L43-L47}{Lean definition} is as follows:
\begin{lstlisting}[language=Lean]
structure PeriodicSpherePacking (d : ℕ) extends SpherePacking d where
  lattice : Submodule ℤ (EuclideanSpace ℝ (Fin d))
  lattice_action : ∀ ⦃x y⦄, x ∈ lattice → y ∈ centers → x + y ∈ centers
  lattice_discrete : DiscreteTopology lattice := by infer_instance
  lattice_isZLattice : IsZLattice ℝ lattice := by infer_instance
\end{lstlisting}
It is easy to see that the density of any sphere packing is at most $1$. Thus, the following supremum, known as the \textbf{sphere packing constant}, is finite:
\begin{align*}
    \Delta_d := \sup_{\Pa(X) \ssq \R^d} \Delta\of{\Pa(X)}
\end{align*}
It is \href{https://github.com/thefundamentaltheor3m/Sphere-Packing-Lean/blob/6bd7a08c1fecdfe6d958574082e8474e602ecf76/SpherePacking/Basic/SpherePacking.lean#L253-L256}{defined in Lean} as follows (where $\bigsqcup$ denotes the \href{https://github.com/leanprover-community/mathlib4/blob/8f9d9cff6bd728b17a24e163c9402775d9e6a365/Mathlib/Order/SetNotation.lean#L58-L61}{(indexed) supremum}):
\begin{lstlisting}[language=Lean]
def SpherePackingConstant (d : ℕ) : ℝ≥0∞ :=
  ⨆ S : SpherePacking d, S.density
\end{lstlisting}
The $E_8$ lattice is the unique even, unimodular lattice in $\R^8$, defined as follows:
\begin{align*}
    \Lambda_8 := \setst{\parenth{x_1, \ldots, x_8} \in \Z^8 \cup \parenth{\Z + \frac{1}{2}}^8}{\sum_{i=1}^{8} x_i \equiv 0 \pmod{2}}
\end{align*}
The $E_8$ sphere packing $\Pa\of{\Lambda_8}$ is the packing of balls of radius $\frac{\sqrt{2}}{2}$ centred at points on $\Lambda_8$. We define it in Lean as the following \lstinline|PeriodicSpherePacking|:
\begin{lstlisting}[language=Lean, escapeinside={(*}{*)}]
noncomputable def E8Packing : PeriodicSpherePacking 8 where
  separation := (* $\sqrt{2}$ *)
  lattice := E8Lattice
  centers := E8Lattice
  centers_dist := by ⋯
\end{lstlisting}
Viazovska's main result in \cite{Viazovska8} is the following:
\begin{theorem}[Viazovska 2016]
    There is no sphere packing in $\R^8$ with density greater than that of the $E_8$ packing. That is,
    \begin{align*}
        \Delta_8 = \Delta\of{\Pa\of{\Lambda_8}} = \pi^4 / 384
    \end{align*}
\end{theorem}
In Lean, our statement of the main theorem is the above equality:
\begin{lstlisting}[language=Lean]
theorem SpherePacking.MainTheorem : SpherePackingConstant 8 = E8Packing.density := ⋯
\end{lstlisting}
The proof relies on the following linear programming bound by Cohn and Elkies. \vspace{-1em}
\begin{theorem}[Cohn-Elkies 2003~{\cite[Theorem 3.1]{CohnElkies}}]\label{Thm:CohnElkies}
    If $f : \R^d \to \R$ is a Schwartz function satisfying the conditions
    \begin{enumerate}[label = \normalfont(CE\arabic*)]
        \item\label{CE1} $f$ is not identically zero.
        \item\label{CE2} For all $x \in \R^d$, if $\norm{x} \geq 1$ then $f(x) \leq 0$.
        \item\label{CE3} For all $x \in \R^d$, $\hat{f}(x) \geq 0$.
    \end{enumerate}
    then we have the following bound on the sphere packing constant $\Delta_d$:
    \begin{align*}
        \Delta_d \leq \frac{f(0)}{\hat{f}(0)} \cdot \Volof{B_d\of{0, \frac{1}{2}}}
    \end{align*}
\end{theorem}
Viazovska \cite[Theorem 3]{Viazovska8} showed the existence of a `magic' function $g : \R^8 \to \R$ whose Cohn-Elkies bound $\frac{g(0)}{\hat{g}(0)} \cdot \Volof{B_d\of{0, \frac{1}{2}}}$ is precisely $\Delta\of{\Pa\of{\Lambda_8}} = {\pi^4}/{384}$.

\subsection{Constructing Viazovska's Magic Function in Lean}

There are deep theoretical challenges to the construction of such a function, and Viazovska's ingenuity was to use quasimodular forms to circumvent them. $g$ is defined in \cite{Viazovska8} as a linear combination of radial, Schwartz Fourier $\pm 1$-eigenfunctions $a$ and $b$ with double zeroes at lattice points, defined as follows:
\begin{align*}
    a(x) :=&
    \int_{-1}^{i} \phi_0\of{\frac{-1}{z + 1}} \, \parenth{z + 1}^2 \, e^{\pi i \norm{x}^2 z} \, \diff{z}
    + \int_{1}^{i} \phi_0\of{\frac{-1}{z - 1}} \, \parenth{z - 1}^2 \, e^{\pi i \norm{x}^2 z} \, \diff{z} \\
    &-2 \int_{0}^{i} \phi_0\of{\frac{-1}{z}} \, z^2 \, e^{\pi i \norm{x}^2 z} \, \diff{z}
    +2 \int_{i}^{i \infty} \phi_0\of{z} \, e^{\pi i \norm{x}^2 z} \, \diff{z}
    \\
    b(x) :=&
    \int_{-1}^{i} \psi_S\of{\frac{-1}{z + 1}} \, \parenth{z + 1}^2 \, e^{\pi i \norm{x}^2 z} \, \diff{z}
    + \int_{1}^{i} \psi_S\of{\frac{-1}{z - 1}} \, \parenth{z - 1}^2 \, e^{\pi i \norm{x}^2 z} \, \diff{z} \\
    &-2 \int_{0}^{i} \psi_S\of{\frac{-1}{z}} \, z^2 \, e^{\pi i \norm{x}^2 z} \, \diff{z}
    -2 \int_{i}^{i \infty} \psi_S\of{z} \, e^{\pi i \norm{x}^2 z} \, \diff{z}
    \\
    g(x) :=& \frac{\pi i}{8640} a(x) + \frac{i}{240\pi} b(x) \qquad \qquad
\end{align*}
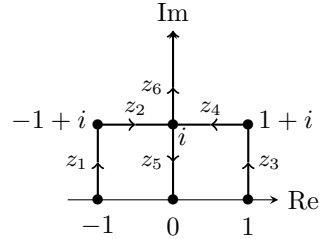
\begin{wrapfigure}[7]{r}{0.36\linewidth}
    \centering
    \vspace{-5.6em}
    \begin{tikzpicture}[scale=1]
        \draw[-stealth] (-1.4,0) -- (1.4,0) node[right] {$\Re$};
        \draw[-stealth] (0,0) -- (0,2.25) node[above] {$\Im$};
    
        \draw[thick, ->] (-1,0) -- (-1,0.5) node[left] {\footnotesize $z_1$};
        \draw[thick, -] (-1,0.5) -- (-1,1);
    
        \draw[thick, ->] (-1,1) -- (-0.5,1) node[above] {\footnotesize $z_2$};
        \draw[thick, -] (-0.5,1) -- (0,1);
        
        \draw[thick, ->] (1,0) -- (1,0.5) node[right] {\footnotesize $z_3$};
        \draw[thick, -] (1,0.5) -- (1,1);
        
        \draw[thick, ->] (1,1) -- (0.5,1) node[above] {\footnotesize $z_4$};
        \draw[thick, -] (0.5,1) -- (0,1);
        
        \draw[thick, ->] (0,1) -- (0,0.5) node[left] {\footnotesize $z_5$};
        \draw[thick, -] (0,0.5) -- (0,0);
    
        \draw[thick, ->] (0, 1) -- (0, 1.5) node[left] {\footnotesize $z_6$};
        \draw[thick, ->] (0, 1.5) -- (0, 2.25);
    
        \labelledpoint{0}{1}{0.125}{-0.6}{$i$}
        \labelledpoint{-1}{0}{0}{-0.8}{$-1$}
        \labelledpoint{0}{0}{0}{-0.8}{$0$}
        \labelledpoint{-1}{1}{-0.65}{-0.4}{$-1 + i$}
        \labelledpoint{1}{0}{0}{-0.8}{$1$}
        \labelledpoint{1}{1}{0.5}{-0.4}{$1 + i$}
    \end{tikzpicture}
    \caption{Choosing contours of integration}
    \label{fig:Contour_param}
\end{wrapfigure}
Denote by $E_2, E_4, E_6$ the weight $2, 4, 6$ Eisenstein Series; $\Delta$ the discriminant form; and $\theta_{10}, \theta_{00}, \theta_{01}$ the Jacobi theta functions. Define $\phi_0$ and $\psi_S$ in terms of these quasimodular forms as follows:
\begin{align*}
    \phi_{0} &:= \frac{\parenth{E_2 E_4 - E_6}^2}{\Delta}
    \\
    \psi_S &:= 128 \parenth{\frac{\theta_{01}^4 - \theta_{10}^4}{\theta_{00}^8} - \frac{\theta_{10}^4 + \theta_{00}^4}{\theta_{01}^8}}
\end{align*}
The above informal definitions of $a$ and $b$ do not specify contours for the constituent integrands. Since the integrands are holomorphic on the upper half-plane, from a strictly mathematical perspective, the choice of contours is immaterial. However, since \mathlib\ lacks sufficient contour integration machinery, we were forced to make a choice of contour when defining $a$ and $b$ in Lean. A crucial step in the proof requires the deformation of unbounded contours (see \cite[\S 4.4]{Sid_Thesis}), and we realized that the best contour shape for such a deformation is rectangular, because \mathlib\ contains a proof of the \href{https://github.com/leanprover-community/mathlib4/blob/19c497800a418208f973be74c9f5c5901aac2f54/Mathlib/Analysis/Complex/CauchyIntegral.lean#L297-L307}{Cauchy-Goursat Theorem for rectangular contours}. We therefore defined rectangular parametrizations $z_1$ to $z_6$ (as shown in \Cref{fig:Contour_param}) and defined $a = I_1 + \cdots + I_6$ and $b = J_1 + \cdots + J_k$, where $I_k$ and $J_k$ are the integrals along $z_k$ with the same integrands as the informal definition. We successfully proved \href{https://github.com/thefundamentaltheor3m/Sphere-Packing-Lean/blob/main/SpherePacking/ForMathlib/CauchyGoursat/OpenRectangular.lean}{unbounded rectangular Cauchy-Goursat}\footnote{Equivalent results were independently proven by contributors to Kontorovich and Tao's open-source \lstinline`PrimeNumberTheoremAnd` formalization project \cite{PNTAnd}}, which was used in \PR{229} to prove that $a$ is expressible in the following form
\begin{align}
    a(x) &= -4 \sinsq{\frac{\pi \norm{x}^2}{2}} \int_{0}^{i\infty} \phi_0\of{\frac{-1}{z}} \, z^2 \, e^{-\pi \norm{x}^2 t} \, \diff{z}
    \label{eq:double_zeroes}
\end{align}

\begin{wrapfigure}[11]{r}{0.36\linewidth}
    \centering
    \vspace{-2.5em}
    \begin{tikzpicture}[scale=1.2]
        \draw[->] (-0.2,0) -- (1.4,0) node[right] {$\Re$};
        \draw[->] (0,-0.2) -- (0,1.4) node[above] {$\Im$};
    
        \draw[thick, domain=0:45, ->] plot ({cos(\x)}, {sin(\x)});
        \draw[thick, domain=45:90, -] plot ({cos(\x)}, {sin(\x)});

        \draw[thick, ->] (1,0) -- (1,0.5);
        \draw[thick, -] (1,0.5) -- (1,1);
        \draw[thick, ->] (1,1) -- (0.5,1);
        \draw[thick, -] (0.5,1) -- (0,1);

        \labelledpoint{0}{0}{-0.3}{-0.8}{$0$}
        \labelledpoint{1}{0}{0}{-0.8}{$1$}
        \labelledpoint{1}{1}{0.3}{-0.2}{$1 + i$}
        \labelledpoint{0}{1}{-0.7}{-0.5}{$i$}
    \end{tikzpicture}
    \caption{Contour needed for eigenfunction property}
    \label{fig:tricky-contour}
\end{wrapfigure}
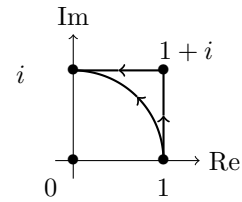

\noindent for all $x \in \R^8$ with $\norm{x} \geq \sqrt{2}$. Replacing $\phi_0$ with $\psi_S$ in \eqref{eq:double_zeroes} yields an analogous expression for $b$. The significance of these representations is that when analytically continued to account for all values of $x \in \R^8$, they yield a convenient expression for $g$ as a certain Laplace transform. A direct computation then shows \ref{CE1}, and the proofs of \ref{CE2} and \ref{CE3} can be reduced, using the facts that $\hat{a} = a$ and $\hat{b} = -b$ (where $\hat{\cdot}$ is the Fourier transform), to the inequalities
\begin{align}
    0 < \phi_0(z) + \frac{36}{\pi^2}\psi_S(z)
    \qquad \qquad \qquad
    0 < \phi_0(z) - \frac{36}{\pi^2}\psi_S(z)
    \label{eqn:modformineqs}
\end{align}
when restricted to the positive imaginary axis \cite[(3), (4)]{Seewoo_Ineq}. The informal and formal machinery for this was largely developed by Lee based on \cite{Seewoo_Ineq}.

Prior to the autoformalization, we did not have a formal proof that $\hat{a} = a$ and $\hat{b} = -b$. These arguments require a version of the Cauchy-Goursat Theorem for contours of the shape depicted in \Cref{fig:tricky-contour} (see \cite[\S 4.3]{Sid_Thesis}), which we had not yet formalized. What made this formalization particularly tricky is that the contour touches the real line, where the integrands, made up of quasimodular forms, were not defined. The saving grace is that in this instance, the integrands vanish because the real line is approached vertically. The Gauss formalization, which shall be discussed in more detail in \Cref{sec:Gauss}, does include a proof that contours of these specific shapes can be deformed. As the project progresses, we would like to generalize the argument to a broader class of contours.

\subsection{The Theory of (Quasi)Modular Forms}

The proof that the magic function satisfies the required properties depends heavily on the basic theory of modular forms. Fortunately, when the Sphere packing project began, \mathlib\ already contained many of the definitions and properties of the Eisenstein series and Jacobi theta functions that we would require. The main development in this area done as part of this project was constructing $E_2$ and $\Delta$ and proving the dimension formulas in level one, which were essential for proving identities between combinations of modular forms, such as $\Delta = (E_4^3 - E_6^2) / 1728$.

The magic function also uses the weight 2 Eisenstein series $E_2$, which is also formalized and already \href{https://github.com/leanprover-community/mathlib4/blob/82df34056bb2de55bb64a01acb804a6a719ef1d0/Mathlib/NumberTheory/ModularForms/EisensteinSeries/E2/Defs.lean#L58-L60}{upstreamed} to \mathlib\ as \href{https://leanprover-community.github.io/mathlib4_docs/Mathlib/NumberTheory/ModularForms/EisensteinSeries/E2/Defs.html#EisensteinSeries.E2}{\lstinline{EisensteinSeries.E2}}. We also formalized the normalized derivative $D := \frac{1}{2\pi i} \frac{\mathrm{d}}{\mathrm{d}z}$ and Serre derivative $\partial_k := D - \frac{k}{12} E_2$. Combined with the dimension formula for level 1 modular forms, we proved Ramanujan's identities on derivatives of Eisenstein series and derivatives of Jacobi theta functions (see \href{https://github.com/thefundamentaltheor3m/Sphere-Packing-Lean/blob/main/SpherePacking/ModularForms/RamanujanIdentities.lean}{\texttt{RamanujanIdentities.lean}} and \href{https://github.com/thefundamentaltheor3m/Sphere-Packing-Lean/blob/main/SpherePacking/ModularForms/ThetaDerivIdentities.lean}{\texttt{ThetaDerivIdentities.lean}} under \href{https://github.com/thefundamentaltheor3m/Sphere-Packing-Lean/tree/main/SpherePacking/ModularForms}{\texttt{SpherePacking/ModularForms}}), which are used to prove various quasimodular form identities used in the proof of the inequalities \cite{Seewoo_Ineq}.
As a by-product, we also proved Jacobi identity $\theta_{00}^4 = \theta_{10}^4 + \theta_{01}^4$.

While the majority of the modular form theory required for this project is standard (see \cite{DiamondShurman} for instance),
a notable obstacle was proving the dimension formula for spaces of modular forms, which, given the lack of contour integration machinery in \mathlib, required a non-standard proof. The details of this and other non-standard proofs will be discussed in a forthcoming paper.

To prove the modular form inequalities \eqref{eqn:modformineqs}, we first defined \lstinline{ResToImagAxis}, the restriction of $F : \mathbb{H} \to \mathbb{C}$ to $\mathbb{R}$, by assigning $F(it)$ for $t > 0$ and 0 otherwise:
\begin{lstlisting}[language=Lean]
def ResToImagAxis (F : ℍ → ℂ) : ℝ → ℂ :=
  fun t => if ht : 0 < t then F ⟨(I * t), by simp [ht]⟩ else 0
\end{lstlisting}
We proved basic properties on realness, positiveness, and how it interacts with (Serre) derivatives (see \projectfile{SpherePacking/ModularForms/ResToImagAxis.lean}). The inequalities \eqref{eqn:modformineqs} will be proved in \PR{331} as \lstinline{FG_inequality_1} and \lstinline{FG_inequality_2}, following the argument of \cite{Seewoo_Ineq}.

\subsection{Metaprogramming}

A rather unexpected outcome of the project so far has been the development of metaprogramming techniques to assist with computations. The first advancement was the development of \lstinline|norm_numI|, a normalization-simplification procedure to bring \lstinline|norm_num|-like functionality to the complex numbers. The tactic recursively parses an expression in $\C$, splits all atoms into their real and imaginary parts, calls \lstinline|norm_num| on their real and imaginary parts, and recombines them using the rules of complex arithmetic to return a normalised, simplified expression in $\C$. The key insight was identifying that the right `normal form' for an expression in $\C$ is $a + bi$, where $a, b \in \R$ are both in normal form (in the sense of \lstinline|norm_num|). The second tactic developed was \lstinline|tendsto_cont|, which uses the continuity of the projection maps from a product of topological spaces to prove goals of the form \lstinline{Tendsto (fun z => expr(f₁ z, ..., fₙ z)) l (nhds c)} where atomic limits \lstinline{Tendsto fᵢ l (nhds aᵢ)} are known from context and \lstinline|expr| is continuous in the atoms (proved via \lstinline{fun_prop}). 
We use it to simplify the proof of the limit of quasimodular forms that are polynomials in $E_2, E_4, E_6, \theta_{10}^4, \theta_{01}^4$.

\section{Autoformalization} \label{sec:Gauss}

Working from the existing blueprint and Lean development, the autoformalization model `Gauss' completed a \lstinline{sorry}-free formalization of the main theorem in five days, taking the project from about $20\,000$ lines to $80\,000$ lines. A subsequent compression and refactoring pass reduced the development to $60\,000$ lines. This archived artifact is available in \PR{341}. Interestingly, it was easier to generate the proof than to clean up and golf the code from 80k lines to 60k. This highlights that the next frontier in autoformalization is not so much getting a proof formalized as doing this while writing high-quality reusable code.

Among the main proofs completed by Gauss were the dimension $8$ contour-deformation identities used in the Fourier eigenfunction argument, including the wedge-set deformation based on the \href{}{Poincar\'e lemma}. In \PR{341}, these appear in \href{https://github.com/math-inc/Sphere-Packing-Lean/blob/gauss/SpherePacking/MagicFunction/a/Eigenfunction/PermI12ContourMain.lean}{\texttt{PermI12ContourMain.lean}} and \href{https://github.com/math-inc/Sphere-Packing-Lean/blob/gauss/SpherePacking/MagicFunction/b/Eigenfunction/PermJ12ContourDeformation.lean}{\texttt{PermJ12ContourDeformation.lean}}. Gauss also completed the finite-dimensionality theory for modular forms on finite-index subgroups, recorded in \href{https://github.com/math-inc/Sphere-Packing-Lean/blob/gauss/SpherePacking/ModularForms/DimensionFormulas.lean}{\texttt{DimensionFormulas.lean}} and \href{https://github.com/math-inc/Sphere-Packing-Lean/blob/gauss/SpherePacking/ModularForms/DimGenCongLevels/FiniteDimensional.lean}{\texttt{FiniteDimensional.lean}}. This general result had already been envisaged in the blueprint as a natural corollary of the surrounding work, although it was not strictly required for the sphere-packing theorem itself.

One stylistic feature of the Gauss contribution is that it introduced many small auxiliary definitions with little or no API attached, many of which a human formalizer would ordinarily inline directly into proofs. A human formalizer would often pause to factor out reusable lemmas and eliminate duplication. Gauss, by contrast, was often content to repeat local patterns. There was a great deal of variation in the mathematical complexity of the Gauss code: it would, at times, write very complex proofs spanning hundreds of lines, yet also prove trivial facts (like $4 = 2 + 2$) or special cases of \mathlib\ results (like Fourier transforms of Gaussians with specific variances). Both extremes are undesirable: long proofs are much harder to maintain and take longer to compile, and reproving results just adds duplication to the codebase. As we clean the Gauss code to bring it up to \mathlib\ standards, we are paying particular attention to these extremes.

\section{A Look to the Future}

As AI models and companies take on more prominent roles in formalization projects, it is worth considering the implications for the practice of mathematics, formal and informal. Some analysis has already been done on the role of Math, Inc in the sphere packing effort \cite{Jeremy_Essay,Simon_Paper} in which both skepticism and optimism are expressed. Of the initial motivations for this project, Avigad writes \cite[p. 3]{Jeremy_Essay},
\begin{quote}
    ``[T]he correctness of Viazovska’s result was never in doubt. The participants embraced the project, rather, as a way of revisiting [Viazovska's] results and better understanding them, and of building libraries and infrastructure to support future work.''
\end{quote}
Gauss, on the other hand, took a more direct route to the verification, which did not completely comply with these objectives. The skepticism on the use of such autoformalization models stems from the initial misalignment of these objectives, and the optimism stems from the belief that the ongoing public collaboration to realign them will motivate the AI for mathematics industry to better consider the objectives of large-scale formalization projects when designing their models.

On a purely technical level, the project was not strictly set up to verify Viazovska's original argument in its entirety, but rather, the modified (if only slightly so) argument in the project blueprint, which differs from the original most notably in its use of Lee's original proof \cite{Seewoo_Ineq} of the inequalities \eqref{eqn:modformineqs}. Furthermore, given the sheer volume of the Gauss code, the exact correspondence between the autoformalization and the intended proof path is still being investigated, and new proofs might still be used to complete certain parts of the project.

Nevertheless, the construction of a \sorry-free proof of the sphere packing problem represents a unique achievement for formal mathematics and machine learning alike. Gauss directly used the human-written Lean code and blueprint, underscoring the importance of robust definitions, detailed proof strategies, and thoughtful design choices in any (auto)formalization workflow. The remainder of the project promises to be immensely interesting, from the formal and informal perspectives alike, and once the project is finished, a complete account of our learnings, along with more technical details of the work done to finish the project, will appear in a forthcoming paper.

\begin{credits}
\subsubsection{\ackname}
The authors are deeply grateful to the many contributors, direct and indirect, to the open-source formalization effort, most notably Matthew Cushman, Cameron Freer, Yongxi Lin, David Loeffler, Heather Macbeth, Pietro Monticone, Aayush Rajasekaran, David Renshaw, Tito Sacchi, Utensil Song, and Yunzhou Xie. The authors also acknowledge the invaluable support of the Lean formalization community, most notably Jeremy Avigad, Kevin Buzzard, Leonardo de Moura, Emily Riehl, and the maintainers of \mathlib. The first author conducted significant portions of this work as an exchange student and subsequently as a research intern at the École Polytechnique Fédérale de Lausanne under the supervision of the seventh author, and later for his master's project \cite{Sid_Thesis} at Imperial College London under the supervision of the fifth author. The first author also acknowledges travel support from the Hoskinson Center for Formal Mathematics at Carnegie Mellon University. The sixth author gratefully acknowledges the support of DARPA's expMath program for the development of Gauss. The sixth and seventh authors acknowledge Renaissance Philanthropy's AI for Math Fund for supporting autoformalization research at EPFL. The authors further acknowledge support from the Institute for Computer-Aided Reasoning in Mathematics on AI policy advice and thank G-Research for sponsoring their weekly `packathons'.

\end{credits}
\bibliographystyle{splncs04}
\bibliography{References}

\end{document}